\documentclass[12pt,twoside,reqno]{amsart}
\usepackage{graphicx}
\usepackage{amsmath}
\usepackage{amsfonts}
\usepackage{amssymb}
\usepackage{amscd}
\usepackage{amsmath}
\usepackage{amscd}
\usepackage{amsthm}
\usepackage{amssymb}
\usepackage{amsfonts}
\usepackage{epsfig}
\usepackage{amsmath,amscd}
\usepackage{graphicx}
\usepackage{geometry,enumerate,latexsym,tabularx,shapepar}
\usepackage{pslatex}
\usepackage[all,2cell,knot,poly]{xy} \UseAllTwocells \SilentMatrices
\usepackage{slashbox}
\input{xypic}
\xyoption{all}
\numberwithin{equation}{section}
\newcommand{\ie}{{\it i.e.\/}\ }

\newcommand{\cf}{{\it cf.\/}\ }

\newtheorem{thm}{Theorem}[section]

\newtheorem{lem}[thm]{Lemma}
\newtheorem{prop}[thm]{Proposition}
\theoremstyle{definition}
\newtheorem{defn}[thm]{Definition}
\theoremstyle{remark}

\numberwithin{equation}{section}
\newtheorem{exa}[thm]{Example}

\newcommand{\bq}{\begin{eqnarray}}
\newcommand{\nq}{\end{eqnarray}}

\newcommand{\m}{\mathbf{M}}
\newcommand{\R}{\mathbb{R}}
\newcommand{\s}{\mathbf{S^3}}

\newcommand{\p}{\pi}

\def\Q{{\mathbb Q}}
\def\Z{{\mathbb Z}}
\def\R{{\mathbb R}}

\def\cK{{\mathcal K}}

\begin{document}

\title [Graph Homologies and Functoriality ]{Graph Homologies and Functoriality}%
\author{Ahmad Zainy Al-Yasry} %
\address{max-planck institute for mathematics \\ and\\ university of baghdad}%
\email{ahmad@azainy.com}%


\date{\today}%


\begin{abstract}
We follow the same technics we used before in \cite{AZ} of extending knot Floer homology
to embedded graphs in a 3-manifold,  by using the Kauffman topological invariant of embedded
graphs by associating family of links and  knots  to a such graph by using some local
replacements at each vertex in the graph.  This new concept of Graph Floer homology
constructed to be the sum of the knot Floer homologies of all the links and knots
associated to this graph and the Euler characteristic is the sum of all the Alexander
polynomials of links in the family. We constructed three pre-additive  categories
one for the graph under the cobordism and the other one is constructed in \cite{AMM}
the last one is a category of Floer homologies for graph defined by Kauffman. Then we
trying to study the functoriality of graphs category and their graph homologies in two
ways, under cobordism and under branched cover, then we try to find the compatibility between them
by the idea of Hilden and Little \cite{HL} by giving a notion of equivalence relation of branched coverings obtained by using
cobordisms, and hence define a a functor from the graph Floer-Kauffman Homology category to the graph Khovanov-Kauffman Homology category.

\end{abstract}

\maketitle

\section{Introduction}

Of course, no introduction can be completely self-contained, so let me say a few words
about the background knowledge the reader is assumed to have.
Many publications appeared studying  Homology in Low dimensional Topology during last two decades.
These studies depended on the concept of categorification, which originally arose in representation theory
and was coined by L. Crane and I. Frenkel \cite{CF}. In the context
of topology, categorification is a process of upgrade algebraic invariants of topological
objects to algebraic categories with richer structure.
The idea of categorification the Jones polynomial is known by
Khovanov Homology for links which is a new link invariant introduced by
Khovanov \cite{kh}, \cite{Ba}, \cite{pt1}. For each link $L$ in $\s$ Khovanov defined a graded chain complex,
with grading preserving differentials, whose graded Euler characteristic is
equal to the Jones polynomial of the link $L$.
In 2003 Ozsv\'{a}th and Szab\'{o} and independently by Rasmussen defined a new useful knot invariant
called  Knot Floer Homology by categorify the Alexander polynomial \cite{os}, \cite{Rs}.
Using grid diagrams for the Heegaard splittings, knot Floer
homology was given a combinatorial construction by Manolescu, Ozsv\'{a}th and Sarkar (2009) \cite{MOS}.
The idea of Khovanov Homology for graphs arises from the same idea of Khovanov homology for links by
the categorifications the chromatic polynomial of graphs.
This was done by L. Helme-Guizon and Y. Rong \cite{laur},
for each graph G, they defined a graded chain complex whose graded
Euler characteristic is equal to the chromatic polynomial of G.
In \cite{AZ} the author defined Khovanov-Kauffman Homology for embedded Graphs by using the
 Kauffman topological invariant of embedded graphs by associating family of links and
 knots to a such graph where the Euler characteristic is the sum of all the Jones polynomials of links in the family.
 Jonathan Hanselman in his work \cite{Jh} described an algorithm to compute $\widehat{HF}$ of any graph
 manifold by using the result of two explicit computations of bordered Heegaard Floer
invariants. The first is the type D trimodule associated to the trivial $S^1$-bundle
over the pair of pants $P$. The second is a bimodule that is necessary for self-gluing,
when two torus boundary components of a bordered manifold are glued to each
other. \\
Yuanyuan Bao in \cite{Yb} defined the Heegaard diagram for a balanced bipartite graph in a rational homology
3-sphere, by introducing a base point for each edge. Then he defined the minus-version and
hat-version of the Heegaard Floer complexes for a given Heegaard
diagram. The hat-version coincides with the sutured Floer complex for the complement
 of the graph, the sutures of which are defined by using the meridians of the edges.
 Bao proved that the homology modules of both versions are topological invariants of the
given graph and discussed some basic properties of the homology. He studied
the Euler characteristic of the hat-version complex. In particular, when the ambient
manifold is the 3-sphere, we give a combinatorial description of the Euler characteristic
by using the "states" of a given graph projection.\\
Shelly Harvey and Danielle O'Donnol \cite{Ho} extend the theory of combinatorial link Floer homology to a class of oriented
spatial graphs called transverse spatial graphs. They defined the notion of a grid diagram
representing a transverse spatial graph, which we call a graph grid diagram. They proved that two
graph grid diagrams representing the same transverse spatial graph are related by a sequence of
graph grid moves, generalizing the work of Cromwell for links. For a graph grid diagram representing
a transverse spatial graph $f : G \rightarrow \s$, they defined a relatively bigraded chain complex (which is a
module over a multivariable polynomial ring) and showed that its homology is preserved under the
graph grid moves; hence it is an invariant of the transverse spatial graph. In fact, they defined both a
minus and hat version. Taking the graded Euler characteristic of the homology of the hat version
gives an Alexander type polynomial for the transverse spatial graph. Specifically, for each transverse
spatial graph $f$, we define a balanced sutured manifold $(\s \setminus f(G); \gamma(f))$. They show that the graded
Euler characteristic is the same as the torsion of $(\s \setminus f(G); \gamma(f))$ defined by S. Friedl, A. Juh\'{a}sz,
and J. Rasmussen.\\
We discuss the question of extending Floer homology from links to
embedded graphs. This is based on a result of Kauffman that constructs a topological
invariant of embedded graphs in the 3-manifold by associating to such a graph
a family of links and knots obtained using some local replacements at
each vertex in the graph. He showed that it is a topological invariant by showing
that the resulting knot and link types in the family thus constructed are invariant
under a set of Reidemeister moves for embedded graphs that determine the
ambient isotopy class of the embedded graphs. We build on this idea and simply
define the Floer homology of an embedded graph to be the sum of the
Floer homologies of all the links and knots in the Kauffman invariant
associated to this graph. Since this family of links and knots is a topologically invariant,
so is the Floer-Kauffman homology of embedded graphs defined in this manner.
\section{Acknowledgements:} The author would like to express his deeply grateful to Prof.Matilde Marcolli
for her advices and numerous fruitful discussions. Some parts of this paper done in Max-Planck Institut f\"ur Mathematik (MPIM), Bonn, Germany,
the author is kindly would like to thank MPIM for their hosting and subsidy during his study there.

\section{ Floer-Kauffman Homology for graphs}
In this part we study a new concept which called Floer-Kauffman Homology for graphs by using Kauffman procedure \cite{kauff} of associating
a family of links to an embedded graph by making some local replacement to each vertex in the graph. We will start recalling
the concept of Floer homology for links (Knots) but not in details and then we combine the definitions of Floer Homology for link with kauffman
idea to produce Floer-Kauffman Homology for graphs. We close this paper by giving an example of computation of Floer-Kauffman homology for an embedded graph using this definition.

\subsection{Link Floer Homology}

 Ozsv\'{a}th - Szab\'{o} and Rasmussen around 2003 Introduced Floer Homology which is an invariant of knots and links in three manifolds. This invariant contain many information about several properties of the knot (genus, slice genus fiberedness, effects of surgery). Ozsv\'{a}th - Szab\'{o}
 used Atiyah Floer Conjecture to develop Heegaard Floer Theory as a symplectic geometric replacement for gauge theory by using Gromov's theory of pseudo-holomorphic curves to construct an invariant of closed 3-manifolds called Heegaard Floer homology. Knot Floer homology is a relative version of Heegaard Floer homology, associated to a pair consisting of a 3-manifolds and a nullhomologous knot in it. Seiberg-Witten equation is the origin of the knot floer homology which is play an important role in 3 and 4 dimensional topology. An invariant called Seiberg-Witten Floer homology of a 3-dimensional manifold $M$ constructed by studying the equation of Seiberg-Witten on the 4-dimensional manifold $M \times \mathbb{R}$. Heegaard Floer homology and Seiberg-Witten Floer homology are isomorphic. Knot Floer homology can be thought of as encoding something about the Seiberg-Witten equations on $\mathbb{R}$ times the knot complement and it is very similar in structure to knot homologies coming from representation
theory, such as those introduced by Khovanov and Khovanov-Rozansky in addition to that it can also be calculated for many small knots using combinatorial methods.\\
Let $K \subset \s$ be an oriented knot. There are several different variants of the knot Floer homology of $K$. The simplest is the hat version, which takes the form of a bi-graded, finitely generated Abelian group

$$\widehat{HFK}(K) = \bigoplus_{ i,s \in \mathbb{Z}} \widehat{HFK}_i(K, s).$$
Here, $i$ is called the Maslov (or homological) grading, and $s$ is called the Alexander grading. The
graded Euler characteristic of $\widehat{HFK}$ is the Alexander-Conway polynomial

$$\sum_{s,i \in\mathbb{Z}} (-1)^iq^s .rank_\mathbb{Z} (\widehat{HFK}_i(K, s)) = \triangle_K(t).$$
\begin{defn}\cite{RA}
The filtered Poincar\'{e} polynomial of a knot homology $\mathcal{H}$ is given by
$$P_{\mathcal{H}(L)}(t, u) = \sum_{i,j}(rank \mathcal{H}_i(L, j)) u^it^j$$
It is a Laurent polynomial in $u$ and $t$.
\end{defn}
If we substitute $u = -1$, the filtered Poincar\'{e} polynomial reduces to the filtered
Euler characteristic. When $rank \mathcal{H}_i(L, j) = 1$, we will often use the shorthand $u^it^j$ to refer to a generator of this group.\\
Let $K$ be a knot in $\s$. There are various Heegaard Floer homology
groups of $K$ such as $\widehat{HFK}$, $HFK^\infty(K)$, $HFK^+(K)$, and $HFK^-(K)$, but our
study will be just with $\widehat{HFK}(K)$. The knot Floer homology $\widehat{HFK}(K)$ is a bigraded chain
complex equipped with a homological grading $u$ and a filtration grading $t$, which
is also known as the Alexander grading. Conventionally, the Alexander grading
is chosen so as to define a downward filtration on $\widehat{HFK}(K)$. The filtered Euler
characteristic of $\widehat{HFK}(K)$ is the Alexander polynomial $\triangle_K(t)$, and the homology
of the complex $\widehat{HFK}(K)$ is a single copy of $\mathbb{Z}$ in homological grading $0$.
We can extend Knot Floer Homology to Links Floer Homology and the Euler characteristic for it multiply by the factor $(t^{1/2} - t^{-1/2})^{n-1}$
where $n$ this the number of the components of the link.
Let $L \subset \s$ be an oriented $n$-component link. Ozsv\'{a}th - Szab\'{o} \cite{os1} show how $L$ can naturally be thought
of as a knot in connected sum of $n-1$ copies of $(S^1 \times S^2)$. This construction gives rise to a knot Floer homology
group $\widehat{HFK}(L)$, which is again a filtered complex. Its filtered Euler characteristic
is given by
$$P_{\widehat{HFK}(L)}(t, -1) = (t^{1/2} - t^{-1/2})^{n-1} \Delta_L(t),$$
and its total homology has rank $2^{n-1}$. The Poincar\'{e} polynomial of the total homology is given by
$$P(u) = (u^{1/2} + u^{-1/2})^{n-1}$$
(when $n$ is odd, the homological grading on $\widehat{HFK}(L)$ is naturally an element of ($\mathbb{Z} + \frac{1}{2}$ rather than of $\mathbb{Z}$.)
In this study, we will not speak more about the many definitions of Knot (Link) Floer homology, but one can
find more details in \cite{os}, \cite{os1}, and \cite{Rs}.
\subsection{Properties}\cite{RA}
Let $L$ be a link in $\s$. Here we give some properties of knot Floer homology.
\begin{prop}\label{1234}
\begin{enumerate}
\item $\widehat{HFK}(L) \cong \widehat{HFK}(L^\circ)$ where $L^\circ$ denotes $L$ with the orientations of all components reversed.
\item $\widehat{HFK}(\overline{L}) \cong \widehat{HFK}(L)^\ast$ where the $\overline{L}$ is the mirror image of $L$, and $\ast$
denotes the operation of taking the dual complex.
\item For two oriented links  $L$ and $L'$, the Knot Floer homology of the disjoint union $ L \sqcup L'$ satisfies
$$\widehat{HFK}(L \sqcup L') = \widehat{HFK}(L) \otimes \widehat{HFK}(L') \otimes X.$$ where $X$ is the rank two complex Poincar\'{e}
$P_X(t,u)= u^{1/2}+u^{-1/2}$ and trivial deferential.
\item  For two oriented links  $L$ and $L'$, the Knot Floer homology of the oriented connected sum  $ L \# L'$ satisfies
$$\widehat{HFK}(L \# L') = \widehat{HFK}(L) \otimes \widehat{HFK}(L').$$
\end{enumerate}
\end{prop}
A natural question is what knot types can be distinguished by knot
Floer homology. If $K_1$ and $K_2$ are distinguished by the Alexander polynomial, then they are also distinguished by knot Floer homology. However, knot Floer homology is a strictly stronger invariant. For example:
\begin{itemize}
  \item If $m(K)$ denotes the mirror of $K$, then $\triangle_K = \triangle_{m(K)}$. On the other hand, $\widehat{HFK}(K) \neq \widehat{HFK}(m(K))$ for the trefoil, and for many other knots;
 \item If $K_1,K_2$ differ from each other by Conway mutation, then $\triangle_{K_1} = \triangle_{K_2}$. A well-known
example of mutant knots, the Conway knot and the Kinoshita-Terasaka knot, have different knot Floer homologies.
\end{itemize}
knot Floer homology is generally an effective invariant for distinguishing between two small knots.
Nevertheless, it has its limitations: we can find examples of different knots with the same knot Floer homology (and, in fact, with the
same full knot Floer complex up to filtered homotopy equivalence).
The alternating knots $7_4$ and $9_2$ are the simplest such example.
A related question is what knots E are distinguished from all other knots by knot Floer homology.
At present, the only known examples are the four simplest knots: the unknot, the two trefoils, and
the figure-eight;
\subsection{Homology theories for embedded graphs}\label{ff}
\begin{figure}
\begin{center}
\includegraphics[width=6cm]{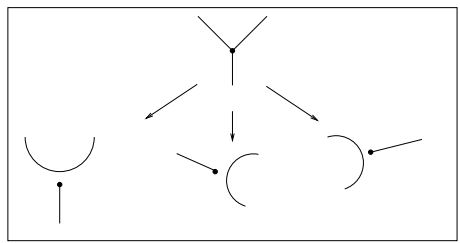}
\caption{local replacement to a vertex in the graph G}\label{fig18}
\end{center}
\end{figure}

\begin{figure}
\begin{center}
\includegraphics[width=4cm]{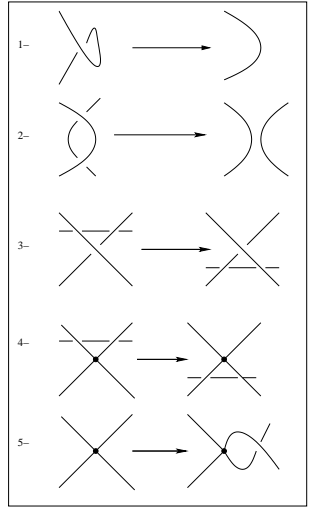}
\caption{Generalized Reidemeister moves by Kauffman}\label{fig21}
\end{center}
\end{figure}

In \cite{kauff} Kauffman introduced an idea relating graphs with links. He by making some local replacements at each vertex in the graph $G$
(Figure\ref{fig18})
can get a family of links associate to graph $G$ which is an invariant under expanded Reidemeister moves defined by Kauffman (Figure \ref{fig21}).
In \cite{AZ} we used kauffman technique to introduce Khovanov-Kauffman homology for embedded graphs. In this part we intend
to define similar homology for graphs by using the concept of Floer homologies for links associated to the graph $G$, to get Floer-Kauffman homology for graphs.

\subsection{Definition of Floer homology for embedded graphs}
We define the concept of Floer homology for embedded graphs by using
Floer homology for the links (knots) and Kauffman theory of associate a family of links
to an embedded graph $G$, as described above.
\subsection{Kauffman's invariant of Graphs}
We give now a survey of the Kauffman theory and show how to associate to an embedded
graph in $3$-Manifold a family of knots and links. We then use these results to give our definition
of Floer homology for embedded graphs.
In \cite{kauff} Kauffman introduced a method for producing topological invariants of graphs
embedded in $3$-Manifold. The idea is to associate a collection of knots and links to a graph $G$ so that
this family is an invariant under the expanded Reidemeister moves defined by Kauffman and
reported here in figure (\ref{fig21}).

He defined in his work
an ambient isotopy for non-rigid (topological) vertices. (Physically, the rigid vertex
concept corresponds to
a network of rigid disks each with (four) flexible tubes or strings emanating from it.)
Kauffman proved that piecewise linear ambient isotopies of embedded graphs in $3$-Manifold
correspond to a sequence of generalized Reidemeister moves for planar diagrams of
the embedded graphs.
\begin{thm} \cite{kauff}
Piecewise linear (PL) ambient isotopy of embedded graphs is generated by the moves of figure (\ref{fig21}), that is, if two embedded graphs
are ambient isotopic, then any two diagrams of them are related by a finite sequence of the moves of figure (\ref{fig21}).
\end{thm}
Let $G$ be a graph embedded in $3$-Manifold. The procedure described by Kauffman of how to
associate to $G$ a family of  knots and links prescribes that we should make a local
replacement as in figure (\ref{fig18}) to each vertex in $G$.
Such a replacement at a vertex $v$ connects two edges and isolates all other edges at that vertex, leaving them as free ends. Let $r(G,v)$ denote the link formed by the closed curves formed by this process at a vertex $v$. One retains the link $r(G,v)$, while eliminating all the remaining unknotted arcs.
Define then $T(G)$ to be the family of the links $r(G,v)$ for all possible replacement choices,
$$ T(G)=\cup_{v\in V(G)} r(G,v). $$

\begin{thm}\cite{kauff}
Let $G$ be any graph embedded in $3$-Manifold, and presented diagrammatically.
Then the family of knots and links $T(G)$,
taken up to ambient isotopy, is a topological invariant of $G$.
\end{thm}

For example, in the figure (\ref{fig1}) the graph $G_2$ is not ambient isotopic to the graph $G_1$,
since $T(G_2)$ contains a non-trivial link.

\begin{defn}\label{401}
Let $G$ be an embedded graph with $ T(G)=\{ L_1,L_2,....,L_n\}$ the family of links associated to $G$
by the Kauffman procedure. Let $\widehat{HFL_r}$ be the usual link Floer homology of the link $L_r$ in this
family. Then the Floer homology
for the embedded graph $G$ is given by $$ \widehat{HFG}=\widehat{HFL_1} \oplus \widehat{HFL_2}\oplus ....\oplus \widehat{HFL_n}$$
Its graded Euler characteristic is the sum of the graded Euler characteristics of the
Floer homology of each link, \ie the sum of the Alexander polynomials,
\begin{equation*}
\sum_{s,i,r \in\mathbb{Z}} (-1)^iq^s .rank_\mathbb{Z} (\widehat{HF_r}_i(K_r, s)) =\sum_r \triangle_{K_r}(q).
\end{equation*}
\end{defn}
We show some simple explicit examples.
\begin{exa}
In figure (\ref{fig1})
\begin{figure}
\begin{center}
\includegraphics[width=8cm]{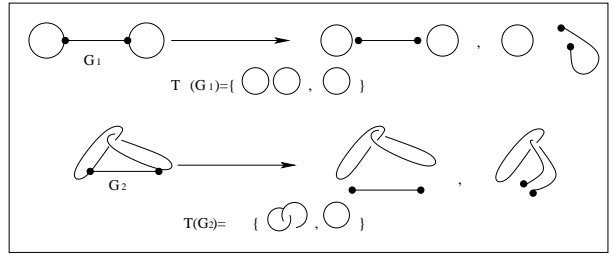}
\caption{Family of links associated to a graph}\label{fig1}
\end{center}
\end{figure}
$T(G_1)=\{\bigcirc \bigcirc ,\bigcirc\}$
$$ \widehat{HFG}(G_1)=\widehat{HFK}(\bigcirc \bigcirc) \oplus \widehat{HFK}(\bigcirc)$$
Now, from proposition \ref{1234} no.3
$$ \widehat{HFG}(G_1)= \widehat{HFK}(\bigcirc) \otimes \widehat{HFK}(\bigcirc) \otimes X \oplus \widehat{HFK}(\bigcirc)$$

Another example comes from $T(G_2)=\{\includegraphics[width=0.6cm]{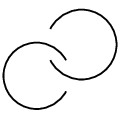}, \bigcirc\}$ then
$$\widehat{HFG}(G_2)=\widehat{HFK}(\includegraphics[width=0.6cm]{hopf.eps}) \oplus \widehat{HFK}(\bigcirc) \otimes X$$

\end{exa}
\newpage
\section{Category of embedded graphs cobordism}
It is interesting to define a cobordism category for graphs by using the idea of kauffman mentioned before. Our construction depends
on associating a family of cobordisms for links to a cobordism for graphs, and hence can define a homomorphism for Khovanov-Kauffman and Floer-Kauffman homologies for graphs.
In the next section by using the idea of kauffman of associating family of links to a graph, we try to define graph cobordism category whose objects are oriented embedded graphs in $\s$, and whose morphisms are isotopy classes of 2-complex oriented, graph cobordisms in
$\s \times [0, 1]$ and hence discuss the functoriality with the category of group of Kauffman-Floer homologies, whose morphisms are isomorphisms.

\subsection{Graph Cobordance}
Here we will define two types of the cobordance between two graphs, one by PL cobordism $\Sigma_P$ defined in (\ref{GraphCobord1}) and the second is by a family of smooth cobordisms $T(\Sigma_\alpha)_{\alpha \in \lambda}$ defined in (\ref{GraphCobord3}). We will prove that these two definitions are
compatible and equivalent by the associativity of the smooth cobordisms $T(\Sigma_\alpha)$ to the PL cobordism $\Sigma_P$. The second definition will help us to define an additive category of embedded Kauffman-graphs under cobordism.
We begin by recalling the definitions of cobordisms in surfaces (smooth and PL surface) and set of smooth surfaces associated to the PL cobordism that will be useful in the rest of our work.
\begin{defn}\label{GraphCobord2}
Two links $ L$ and $L'$ are called cobordic if there is a surface $\Sigma$ have the
boundary $\partial \Sigma= L \cup -L'$ with $L= \Sigma \cap (\s\times \{0\})$, $L'=\Sigma \cap (\s\times \{1\})$. Here by "surfaces" we mean
$2$-dimensional compact differentiable manifold embedded in $\s \times [0,1]$. We define the identity cobordism to be $id_L$ for a link $L$. {\em $[\Sigma_{L}]$ denotes the cobordism class of the link $L$.}
\end{defn}
\begin{defn}\label{GraphCobord1}\cite{AZY}
Two graphs $ G$ and $G'$ are called cobordic if there is a PL surface $\Sigma_P$ have the
boundary $\partial \Sigma_P= G \cup -G'$ with $G= \Sigma \cap (\s\times \{0\})$, $G'=\Sigma \cap (\s\times \{1\})$. Here by "surfaces" we mean
$2$-dimensional simplicial complexes that are PL-embedded in $\s \times [0,1]$. We define the identity cobordism to be $id_G$ for a graph $G$. {\em $[\Sigma_{P}]$ denotes the cobordism class of the graph $G$.}
\end{defn}
\begin{defn}\label{GraphCobord3}
 Let $G, G'$ two embedded graphs have associated families of links $T(G)$ and $T(G')$ respectively according to Kauffman construction, and let $T(\Sigma_\alpha)_{\alpha \in \lambda}$  be a family of smooth cobordisms have links boundaries of $T(G)$ and $T(G')$.
These graphs are said to be cobordant if there is a family of smooth cobordisms $T(\Sigma_\alpha)_{\alpha \in \lambda}$ in $\s \times [0,1]$ such that $T(G)$ and $T(G')$ are the boundary of $T(\Sigma_\alpha)$. These family of smooth cobordisms defined in (\ref{GraphCobord2}) can be associated to the PL cobordism $\Sigma_P$ between $G$ and $G'$ defined in (\ref{GraphCobord1}).
\end{defn}
\begin{thm}\label{assotoplcobo}
The family of smooth cobordisms $T(\Sigma_\alpha)_{\alpha \in \lambda}$ is associated to the PL cobordism $\Sigma_P$.
\end{thm}
\proof
We show that the union $L:=\cup_{t\in [0,1]} L_t$ is a cobordism between the Kauffman invariants $L_0$ of the graph $G_0$ and $L_1$ of the graph $G_1$. Let $\Sigma_P$ be a 2-complex embedded in $\s \times [0,1]$, with boundaries $G_0$ and $G_1$. Define the map $p:\s \times [0,1] \rightarrow [0,1]$ to be the projection on the second factor. Then, for each $t \in [0,1]$, consider $G_t:=p^{-1}(t)\cap \Sigma_P$. This is an embedded graph in $\s$, You can apply to $G_t$ the construction of Kauffman to get $L_t$ (a collection of links associated to $G_t$ for each $t$).
 By using Kauffman procedure we can associate for all $G_t$ and  $t \in [0,1]$ a family of links and this means we can obtain many cobordisms between $L_0$ and $L_1$ in different ways but they are all equivalent and hence we can make a class of cobordisms for links $[Cob L]$ between $L_0$ and $L_1$. In case the 2-complex for some $t \in [0,1]$ be a point or line or even a graph but has no family of links associate this means that the set is empty $\phi$ and for some $(\epsilon  \in R)$ where $G_{t - \epsilon}$ and $G_{t + \epsilon}$  graphs have families of links then we can use unit and counit to close the smooth surface.
 \endproof
\subsection{Composition of Cobordisms}
We can define the composition between two PL cobordisms $\Sigma_P$ and $\Sigma'_P$ (where $\Sigma_P$ is the cobordism with boundary
$G \cup -G'$,  $\Sigma'_P$ is the cobordism with boundary $G'\cup -G''$ in $\s \times [0,1]$) to be $\Sigma_P \circ \Sigma'_P$ which is a PL cobordism with boundary $G \cup -G''$.
For the second type of the cobordance by family of smooth cobordisms we can define the composition as follows :
Let $G, G'$ and $G''$ be embedded graphs in $\s$ with $T(G), T(G')$ and $T(G'')$ links families associated to each graph respectively.
We have a family of smooth cobordisms $T(\Sigma_\alpha)_{\alpha \in \lambda}$ have boundary $T(G) \cup -T(G')$ and another family
$T(\Sigma'_\beta)_{\beta \in \gamma}$ with boundary $T(G') \cup -T(G'')$. Let $\Sigma_\beta$ and $\alpha \in \lambda$ be a smooth cobordism
with boundary links from the sets $T(G)$ and $T(G')$ and let $\Sigma'_\beta$ for $\beta \in \gamma$ be a smooth cobordism
with boundary links from the sets $T(G')$ and $T(G'')$. $\Sigma_\alpha \circ \Sigma'_\beta$ is a smooth cobordism with boundary links
from the sets $T(G)$ and $T(G'')$ and this defines the composition of the second type of the graphs cobordance.
Let $\Sigma_P \in Hom (G,G')$ be a graph cobordism, we can think of the class of link cobordisms as a morphism between $G$ and $G'$
$$\Sigma_P \equiv \sqcup_{\alpha \in \lambda} \Sigma_\alpha$$
$\Sigma_P$ is not unique, we can get another $\widetilde{\Sigma_P} \in Hom (G,G') $  with the same
class of link cobordism $[\Sigma_{T(G)}]$
with the addition of coefficients.
Let $\mathcal{G_K}$ be a category, whose objects are embedded graphs that have a family of links according to the kauffman definition and
morphisms are $2$-dimensional simplicial complex surface $\Sigma_P \in Hom (G,G')$ with graphs boundary. As we maintained before, to each graph cobordism we can associate a family of link cobordisms with boundary family of links associated to each graph. There are many link cobordism can be associated but they are all equivalent.

\begin{defn}\label{preadd}
A pre-additive category $\mathcal{C}$ is a category such that, for
any $\mathcal{O},\mathcal{O'} \in Obj(\mathcal{C})$ the set of
morphisms $Hom(\mathcal{O},\mathcal{O'})$ is an abelian group and
the composition of maps is a bilinear operation, that is, for
$\mathcal{O},\mathcal{O'}, \mathcal{O''} \in Obj(\mathcal{C})$ the
composition
$$ \circ :
Hom(\mathcal{O},\mathcal{O'})\otimes Hom(\mathcal{O'},\mathcal{O''}) \longrightarrow
Hom(\mathcal{O},\mathcal{O''})$$ is a bilinear homomorphism.
\end{defn}
\begin{lem}
$\mathcal{G_K}$ is a pre-additive category.
\end{lem}
In fact, we can write PL morphisms $\Sigma_P$ between graphs in term of smooth morphisms $\Sigma_P$ between links as follows:
$$\Sigma_P=\sum_i a_i (\Sigma_P)_i$$ equivalently as $$\Sigma_P=\sum a_n \Sigma_P$$ where
the sum ranges over the set of all PL cobordisms and
all but finitely many of the coefficients $a$ are zero. Then, for
$$\Sigma_P=\sum a_n \Sigma_P$$ and $$\widetilde{\Sigma_P}=\sum b_n \Sigma_P$$ we have
$$\Sigma_P+\widetilde{\Sigma_P}=\sum (a_n + b_n) \Sigma_P$$
The composition rule given by composition of cobordisms. We can compose cobordisms by gluing graphs. This gives a bilinear
homomorphism
$$ Hom(G,G')\otimes Hom(G',G'')\to Hom(G,G''). $$
This shows that $\mathcal{G_K}$ is a pre-additive category.
\endproof
In our case, the set of morphisms $Hom(G,G')$ is an abelian group
\begin{defn}
Suppose given a pre-additive category $\mathcal{C}$. Then the
additive category $Mat(\mathcal{C})$ is defined as follows (\cf
\cite{Ba}).
\begin{enumerate}
\item The objects in $Obj(Mat(\mathcal{C}))$ are formal direct
sums $\bigoplus_{i=1}^{n} \mathcal{O}_i$ of objects $\mathcal{O}_i
\in Obj(\mathcal{C})$, where we allow for the direct sum to be
possibly empty.
\item If $F:\mathcal{O'}\rightarrow \mathcal{O}$ is a morphism in $Mat(
\mathcal{C})$ with objects $\mathcal{O}=\bigoplus_{i=1}^m
\mathcal{O}_i$ and $\mathcal{O'}=\bigoplus_{j=1}^n \mathcal{O}_j$
then $F=F_{ij}$ is a $m\times n$ matrix of morphisms
$F_{ij}:\mathcal{O'}_j\rightarrow \mathcal{O}_i$ in $\mathcal{C}$.
The abelian group structure on
$Hom_{Mat(\mathcal{C})}(\mathcal{O'},\mathcal{O})$ is given by
matrix addition and the abelian group structure of
$Hom_{\mathcal{C}}(\mathcal{O'}_j,\mathcal{O}_i)$.
\item The composition of morphisms in $Mat(\mathcal{C})$ is defined by
the rule of matrix multiplication and the composition of morphisms
in $\mathcal{C}$.
\end{enumerate}
Then $Mat(\mathcal{C})$ is called the {\em additive closure} of
$\mathcal{C}$. For more details see for instance \cite{Ba}.
\end{defn}
In the following, for simplicity of notation, we continue to use the
notation $\mathcal{G_K}$ for the additive closure of the category $\mathcal{G_K}$ of
Definition \ref{catknot1}.
Define an equivalence relation between two graphs $G$ and $G'$ as follows: two graphs are said to be equivalent $G \sim G'$ if they have
the same set of links associated to both $G$ and $G'$ for each local replacement to each vertex in $G$ and $G'$, and hence their Floer-Kauffman homologies are isomorphic $\widehat{HFG}(G) \simeq \widehat{HFG}(G')$. $[G]_K$ is the equivalence class
of graphs under the Kauffman idea. This equivalence built a new class for Khovanov-Kauffman and Floer-Kauffman homologies.
These definitions of cobordance of $G$ and $G'$ can define an equivalence classes for graphs. Two graphs $G$ and $G'$ are said to be equivalent if there is $\Sigma_P$ PL cobordism (or family of smooth cobordisms between ) with boundary $G \cup -G'$. $G$ and $G'$ are equivalent (and have the same set of links $T(G) \equiv T(G')$) then graph Khovanov-Kauffman and Floer-Kauffman homologies are equivalent (\ie $FKh(G)\equiv FKh(G')$ or $KKh(G)\equiv KKh(G')$), and this class induce a Khovanov morphism $Kh(L_0) \to Kh(L_1)$ also this class of cobordisms can be associate to the 2-complex with boundary $G_0$ and $G_1$.
\subsection{Category of embedded graphs under branched covers}
In \cite{AMM} we constructed a Statistical Mechanical System of 3-manifolds as branched
coverings of the 3-sphere, branched along embedded graphs (or in particular knots) in the 3-sphere
and where morphisms are formal linear combinations of 3-manifolds. Our definition of correspondences
reliesed on the Alexander branched covering theorem, which shows that all compact
oriented 3-manifolds can be realized as branched coverings of the 3-sphere, with branched locus an
embedded (not necessarily connected) graph. The way in which a given 3-manifold is realized as
a branched cover is highly not unique. It is precisely this lack of uniqueness that makes it possible
to regard 3-manifolds as correspondences.
\begin{defn}\cite{pr}\label{branch}
A branched covering of 3-manifolds is defined as a
continuous map $p : M^3 \to  N^3$ such that there exists a one-dimensional subcomplex
$l^1$ in $N^3$ whose inverse image $p^{-1}(l^1)$ is a one-dimensional subcomplex on the
complement to which, $M^3 - p^{-1}(l^1)$, the restriction of $p$ is a covering. In this
situation $M^3$ is called the covering manifold, $N^3$ is the base, and $l^1$ is the branching
set.
\end{defn}
Using cyclic branched coverings $p : \s \to \s$ as the starting point, it is
possible to construct other examples of branched coverings by performing surgery
along framed links. In the base of the branched covering $p$, let us choose
a framed link $L$ and do surgery along $L$, producing a manifold $N^3$. Consider the
inverse image of $L$ under the map $p$ and perform surgery along it, obtaining another
manifold $M^3$. The branched covering $p : \s \to \s$ induces the branched covering
$\pi : M^3 \to N^3$; we shall also call such a branched covering cyclic.
\begin{exa}\label{Psphere} {\em (Poincar\'e homology sphere):} Let
$\mathbf{P}$ denote the Poncar\'e homology sphere. This smooth
compact oriented 3-manifold is a 5-fold cover of $\s$ branched along
the {\em trefoil knot} (that is, the $(2,3)$ torus knot), or a
3-fold cover of $\s$ branched along the $(2,5)$ torus knot, or also
a 2-fold cover of $\s$ branched along the $(3,5)$ torus knot. For
details see \cite{pr}, \cite{ks}.
\end{exa}
We constructed an additive category whose objects are embedded graphs in the
3-sphere and where morphisms are formal linear combinations of 3-manifolds.
In this section we will recall the definitions of morphisms between graphs and trying
to establish later a functor with the graph homologies by using the idea of  Hilden and Little (\cite{HL}).\\
 Define $\phi:G \rightarrow G'$ between
graphs as formal finite linear combinations
\bq \label{phiMi}
\phi=\sum_{i} a_{i}  \mathbf{M}_i
\nq
with $a_{i} \in \Q$ and
$\mathbf{M}_i$ compact oriented smooth 3-manifolds with submersions
$$ \pi_{i}:\m_{i} \rightarrow \s $$ and $$ \pi'_{i}:\m_{i}
\rightarrow \s $$ that are branched covers, respectively branched
along $G$ and $G'$.
\begin{exa}\label{cyclexa} {\em (Cyclic branched coverings):}
We represent $\s$ as $\R^3 \cup \{ \infty \}$. Let $l$ be a straight
line chosen in $\R^3$. Consider the quotient map $p:\R^3 \rightarrow
\R^3/(\Z/n\Z)$ that identifies the points of $\R^3$ obtained from
each other by a rotation by an angle of $\frac{2\pi}{n}$ about the
axis $l$. Upon identifying $\R^3 \simeq \R^3/(\Z/n\Z)$, this extends
to a map $p:\s \rightarrow \s$ which is an $n$-fold covering
branched along the unknot $l\cup \{ \infty \}$ and with multiplicity
one over the branch locus.
\end{exa}
We use the notation \bq \label{2covers} G \subset \s
\stackrel{\pi_G} \longleftarrow \mathbf{M} \stackrel{\pi_{G'}}
\longrightarrow \s \supset G' \nq for a 3-manifold that is realized
in two ways as a covering of $\mathbf{S^{3}}$, branched along the
graph $G$ or $ G'$.
These cyclic branched coverings are useful to construct other more
complicated branched coverings by performing surgeries along framed
links (\cite{pr}).
\begin{defn}\cite{AMM}\label{fibproddef}
Suppose given \bq  \label{MMprime} G \subset \s \stackrel{\pi_G}
\longleftarrow \m \stackrel{\pi_{G'}} \longrightarrow \s \supset G'
\ \ \ \text{ and }  \ \ \ G'\subset \s \stackrel{\tilde{\pi}_{G'}}
\longleftarrow \tilde\m \stackrel{\tilde{\pi}_{G''}} \longrightarrow
\s \supset G'' . \nq One defines the composition $\m \circ \tilde\m$
as \bq \label{McircMprime} \m \circ \tilde{\m} :=\m \times_{G'} \tilde\m ,
\nq where the fibered product $\m \times_{G'} \tilde\m$ is defined
as \bq \label{fibprod} \m \times_{G'} \tilde\m :=\{(x,x')\in \m
\times \tilde\m |\pi_{G'}(x)=\tilde\pi_{G'}(x')\}. \nq
\end{defn}
The composition $\m \circ \tilde\m$ defined in this way satisfies
the following property.

\begin{prop} \cite{AMM}\label{compbranch}
Assume that the maps of \eqref{MMprime} have the following
multiplicities. The map $\p_G$ is of order $m$ for $x\in \s
\smallsetminus G$ and of order $n$ for $x\in G$; the map $\p_{G'}$
is of order $m'$ for $x\in \s \smallsetminus G'$ and $n'$ for $x\in
G'$; the map $\tilde{\p}_{G'}$ is of order $\tilde m'$ for $x\in \s
\smallsetminus G'$ and of order $\tilde n'$ for $x\in G'$; the map
$\tilde{\p}_{G''}$ is of order $\tilde m''$ for $x \in \s
\smallsetminus G''$ and $\tilde n''$ for $x\in G''$.
For simplicity assume that  \bq \label{emptyint} G \cap
\pi_G(\pi_{G'}^{-1}(G')) =\emptyset \ \ \ \text{ and }  \ \ \ G''
\cap \tilde\pi_{G''}(\tilde\pi_{G'}^{-1}(G'))=\emptyset. \nq
Then the fibered product $\hat{\m}=\m \times_{G'} \tilde\m$ is a
smooth 3-manifold with submersions \bq \label{coverhatM} E \subset
\s \stackrel{\hat{\pi}_E} \longleftarrow \hat{\m}
\stackrel{\hat{\pi}_{E''}} \longrightarrow \s \supset E'' . \nq
where
\bq E=G \cup  \pi_{G}(\pi_{G'}^{-1}(G'))\\
      E''=  G'' \cup  \tilde \pi_{G}(\tilde\pi_{G'}^{-1}(G'))
\nq
\end{prop}
This definition makes sense, since the way in which a given
3-manifold $\mathbf{M}$ is realized as a branched cover of $\s$
branched along a knot is not unique.
\begin{defn}\cite{AMM} \label{catknot1}
We let $\cK$ denote the category whose objects $Obj(\cK)$ are graphs
$G\subset \s$ and whose morphisms $\phi\in Hom(G,G')$ are
$\mathbb{Q}$-linear combinations $\sum_{i} a_i \m_i$ of 3-manifold
$\m_i$ with submersions $\p_{E}$ and $\p_{E'}$ to $\s$, including the trivial (unbranched)
covering in all the $Hom(G,G)$.
\end{defn}

\begin{lem}\cite{AMM}\label{preadditive}
The category $\cK$ is a small pre-additive category.
\end{lem}
The proof is in \cite{AMM}.

\subsection{Category of graph Khovanov and Floer homologies}
In this section we define 4 homology groups categories (two for links homologies and two graphs homologies). Let $\mathcal{C}_{Kh}$ be a category of Khovanov homology for links whose objects are khovanov homology groups and morphisms are group homomorphisms $f_{Kh}: Kh(L) \to Kh(L')$ for some $L, L'$ links. In a parallel direction we can define a category $\mathcal{C}_{Fh} $ of Floer homologies for links, whose objects are links floer homology group defined in section (\ref{ff}) and morphisms are group homomorphisms $f_{Fh}: Kh(L) \to Kh(L')$ for some $L, L'$ links. Clearly $\mathcal{C}_{Kh}$ and $\mathcal{C}_{Fh}$ are Additive categories.
In \cite{AZ} we introduced the notation of Khovanov-Kauffman homology for graphs. To each embedded graph who can be under a certain of local replacements at each vertex associated by a family of links or knots, this definition help us to define a new category of Khovanov-Kauffman homology for groups $\mathcal{C}_{KKh}$ whose objects are Khovanov-Kauffman homology groups and morphisms are group homomorphisms $f_{KKh}: KKh(G) \to KKh(G')$ for some $G, G'$ graphs. Here $K$ refereing to the kauffman. In definition (\ref{401}) we defined Floer-Kauffman homology for graphs and hence we can define and new category of Khovanov-Kauffman homology for groups $\mathcal{C}_{FKh}$ whose objects are Khovanov-Kauffman homology groups and morphisms are group homomorphisms $f_{FKh}: FKh(G) \to FKh(G')$ for some $G, G'$ graphs.

\begin{lem}\label{40000}
The categories $\mathcal{C}_{KKh}$ and $\mathcal{C}_{FKh}$ are small pre-additive (Additive) categories.
\end{lem}
In the next section we need to study the Functoriality between the cobordisms categories and the homology categories.


\section{Functoriality}
\subsection{Functoriality under PL Cobordisms}
In this part we want to study the Functoriality between link, graph categories in this side and their homology groups categories in the other side.
\subsubsection{Functoriality under Khovanov homology}\label{KKH}
An important problem in this theory is to extend the Khovanov homology to a
monoidal functor from the category of cobordisms of oriented links. The first attempt
by Khovanov gave a negative answer because of many problems of signs. Functoriality
up to sign was conjectured by Khovanov and proved later by Jacobson \cite{Ja}, Bar Natan
\cite{BN2} and Khovanov \cite{Kh3}. This functoriality up to sign was used by Rasmussen \cite{jr} to prove a conjecture of Milnor about the slice genus.
Let $[\mathcal{C}_L]$ denoted the link cobordism category whose objects are oriented links in $\s$, and whose morphisms are isotopy classes of oriented, link cobordisms in $\s \times [0, 1]$. In this section we trying to introduce the functoriality between graphs cobordism category and Khovanov-Kauffman homology under the concept of PL cobordism. We give now the definition of cobordism class of the graph $G$. We can define the morphism of Khovanov-Kauffman homology for graphs $KKh(G_0) \to KKh(G_1)$ to induce the functoriality by using the same procedure prove the functoriality of Floer-Kauffman homologies obtain a morphism between the homologies using the cobordism L and existing results on
 functoriality of Khovanov or Floer theory of links with respect to link cobordisms.
Let $G$ and $G'$ be two embedded graphs in $\s$ and $\Sigma_P$ is a PL cobordism in $\s \times [0, 1]$ with $\partial \Sigma_P = G \cup -G'$.
In \cite{Jaco} Jacobsson showed that a movie for a cobordism in $\s \times [0, 1]$ with starting and ending diagrams $D_0$ and $D_1$ for links $L_0$ and $L_1$ respectively can induces a map

\begin{equation}\label{1000}
  f_D: Kh(D_0) \to Kh(D_1)
\end{equation}
Where $Kh(.)$ is the khovanov homology. Two movies from $D_0$ to $D_1$ supposed to be equivalent if their lifts, for fixed links $L_0$ and $L_1$, represent the same morphism in $[\mathcal{C}_L]$, and this means Khovanov homology is really a functor
\begin{equation}\label{1111}
    Kh : \mathcal{C}_L \to \mathcal{V}
\end{equation}

where $[\mathcal{V}]$ is the category of Vector spaces.\\

If we have a family of smooth cobordism $T(\Sigma_\alpha)_{\alpha \in \lambda}$ then we can get a family of morphisms $T(f_\alpha)_{\alpha \in \lambda}$ where $f: Kh(L_0) \to Kh(L_1)$ is a homomorphism between khovanov homologies for links $L_0$ and $L_1$.
Since $G$ and $G'$ are two embedded graphs in $\s$, then according to kauffman construction we can associate
a family of links to both $G$ and $G'$ and denoted $T(G), T(G')$. We defined Khovanov-Kauffman homology for graphs as:
$KKh(G)= \oplus^n_1 Kh(L_i)$ for $1 \leq i \leq n$ and $KKh(G')= \oplus^m_1 Kh(L'_j)$ for $1 \leq j \leq m$ where $n,m \in \mathbb{Z}$.
In theorem (\ref{assotoplcobo}) we showed that, for a PL cobordism $\Sigma_P$ we can associate a family of smooth cobordism $T(\Sigma_\alpha)_{\alpha \in \lambda}$ to $\Sigma_P$ and hence a family of khovanov homology homomorphisms induced by these smooth cobordisms, this family of homomorphisms induced by smooth cobordisms will induce a homomorphism $f_P : KKh(G) \to KKh(G')$ from a PL cobordism $\Sigma_P$ with boundary $G,G'$.
        \[
\begin{xy}
\xymatrix{
 T(\Sigma_\alpha)_{\alpha \in \lambda} \ar[d]^{induce} \ar[r]^{asso} &  \ar[d]^{induce}\Sigma_P \\
T(f_\alpha)_{\alpha \in \lambda} \ar[r]^{asso} & f_P }
\end{xy}
\]
For some interval $[t,t+\varepsilon] \subseteq [0,1]$ where $\varepsilon \in \mathbb{R}$ let $\Sigma_{t,P}$ be a partial PL cobordism in $\Sigma_{P}$ with boundary graphs $\partial\Sigma_{t,P} = G_t \cup -G_{t+\varepsilon} $. If both $G_t$ and $G_{t+\varepsilon} $ have family of links then we can
define a family of homomorphisms as we explained before. If $G_t$ and $G_{t+\varepsilon} $ have no or one of them have and the other not this means
we have unit and counit cobordisms then where they define homomorphisms to complex field $\mathbb{C}$.
\subsubsection{Functoriality under Floer homology}\label{FKH}
It has been a central open problem in Heegaard Floer theory whether cobordisms of links induce homomorphisms on the associated link
Floer homology groups. In \cite{AJU}  Andr\'{a}s Juh\'{a}sz established a cobordism between sutured manifolds and showed that such
a cobordism induces a map on sutured Floer homology. This map is a common
generalization of the hat version of the closed 3-manifold cobordism map in
Heegaard Floer theory, and the contact gluing map. Since the introduction of knot Floer homology, it has been a natural question
whether knot cobordisms induce maps on knot Floer homology, exhibiting it as
a categorification of the Alexander polynomial.  We show that cobordisms of sutured manifolds induce maps on sutured
Floer homology, a Heegaard Floer type invariant of 3-manifolds with boundary. Cobordism maps in Heegaard Floer homology were
first outlined by Ozsv\'{a}th and Szab\'{o}  for cobordisms between closed 3-manifolds,
but their work did not address two fundamental questions. The first was the issue
of assigning a well-defined Heegaard Floer group – not just an isomorphism class
– to a 3-manifold, and the functoriality of the construction under diffeomorphisms. The second issue was exhibiting the
independence of their cobordism maps of the surgery description of the cobordism.
They did check invariance under Kirby moves, but did not address how this gives
rise to a well-defined map without running into naturality issues. Sutured manifolds, introduced by Gabai \cite{GAA}, have been of great use in 3-manifold
topology, and especially in knot theory. A sutured manifold $(M, \gamma)$ is a compact
oriented 3-manifold M with boundary, together with a decomposition of the boundary
$\partial M$ into a positive part $R_+(\gamma)$ and a negative part $R_-(\gamma)$ that meet along a
"thickened" oriented 1-manifold $\gamma \subset \partial M$ called the suture. A. Juh\'{a}sz \cite{HHH}, \cite{HHT} defined an invariant called sutured Floer homology, in short $SFH$, for balanced sutured manifolds. $SFH$ can be viewed as a common generalization of
the hat version of Heegaard Floer homology and link Floer homology, both defined
by Ozsv\`{a}th and Szab\'{o}. A. Juh\'{a}sz showed that $SFH$ behaves nicely
under sutured manifold decompositions, which has several important consequences,
such as the above-mentioned detection of the genus and fibredness by knot Floer
homology.
Ozsv\'{a}th - Szab\'{o} in \cite{PPZ} introduced Heegaard Floer homology of branched double cover. For any $p \in \s$ They defined a
Based link cobordism category $\mathcal{L}_p$ where objects are oriented based links in $\s$ containing point $p$ and
whose morphisms are isotopy classes of oriented link cobordisms in $\s \times [0,1]$. Given a based link cobordism $S$ from $L_0$
to $L_1$, the branched cover $\s \times [0,1]$ along $S$ is a smooth oriented 4-dimensional cobordism $\Sigma(S)$ from $\Sigma(L_0)$ to
$\Sigma(L_1)$ and induce a map on Heegaard Floer Homology
$$ \widehat{HF}(-\Sigma(S)): \widehat{HF}(-\Sigma(L_0)) \to \widehat{HF}(-\Sigma(L_1)) $$
which is an invariant of the morphism in $\mathcal{L}_p$ represented by $S$. The Heegaard Floer homology of branched
cover can define a functor $$\mathcal{F}: \mathcal{L}_p \to \mathcal{V}$$ where $\mathcal{V}$ is the category of vector spaces.

If we have a family of smooth cobordism $T(\Sigma_\alpha)_{\alpha \in \lambda}$ then we can get a family of morphisms $T(f_\alpha)_{\alpha \in \lambda}$ where $\widehat{HFG}(G_1): \widehat{HFK}(L_0) \to \widehat{HFK}(L_1)$ is a homomorphism between Floer homologies for links $L_0$ and $L_1$. Since $G$ and $G'$ be two embedded graphs in $\s$ then according to kauffman construction we can associate
a family of links to both $G$ and $G'$ and denoted $T(G), T(G')$. We define the Floer-Kauffman homology for graphs as:
$\widehat{HFG}(G)= \oplus^n_1 \widehat{HFK}(L_i)$ for $1 \leq i \leq n$ and $\widehat{HFG}(G')= \oplus^m_1 \widehat{HFK}(L'_j)$ for $1 \leq j \leq m$ where $n,m \in \mathbb{Z}$.
In theorem (\ref{assotoplcobo}) we showed that, for a PL cobordism $\Sigma_P$ we can associate a family of smooth cobordism $T(\Sigma_\alpha)_{\alpha \in \lambda}$ to $\Sigma_P$ and hence a family of Floer-Kauffman homomorphisms induced by these smooth cobordisms this family of homomorphisms induced
by smooth cobordisms will induce a homomorphism $f_P : \widehat{HFG}(G) \to \widehat{HFG}(G')$ from a PL cobordism $\Sigma_P$ with boundary $G,G'$.
        \[
\begin{xy}
\xymatrix{
 T(\Sigma_\alpha)_{\alpha \in \lambda} \ar[d]^{induce} \ar[r]^{asso} &  \ar[d]^{induce}\Sigma_P \\
T(f_\alpha)_{\alpha \in \lambda} \ar[r]^{asso} & f_P }
\end{xy}
\]
For some interval $[t,t+\varepsilon] \subseteq [0,1]$ where $\varepsilon \in \mathbb{R}$ let $\Sigma_{t,P}$ be a partial PL cobordism in $\Sigma_{P}$ with boundary graphs $\partial\Sigma_{t,P} = G_t \cup -G_{t+\varepsilon} $. If both $G_t$ and $G_{t+\varepsilon} $ have family of links then we can
define a family of homomorphisms as we explained before. If $G_t$ and $G_{t+\varepsilon} $ have no or one of them have and the other not this means
we have unit and counit cobordisms then where they define homomorphisms to complex field $\mathbb{C}$.

\section{Compatibility between functoriality under cobordism and branched cover}\label{HL}
In \cite{AMM} we constructed cobordism of branched cover by using the work of Hilden and Little (\cite{HL}).
They gave us a suitable notion of equivalence relation of branched coverings obtained by using
cobordisms. In this part we need to invest this notation to study the compatibility between the
functoriality under the cobordism and branched cover.
Suppose given two compact oriented 3-manifolds
$\m_1$ and $\m_2$ that are branched covers of $\s$, with covering
maps $\pi_1: \m_1 \to \s$ and $\pi_2: \m_2 \to \s$,  respectively
branched along $1$-dimensional simplicial complex  $E_1$ and $E_2$. A cobordism of branched
coverings is a $4$-dimensional manifold $W$ with boundary $\partial
W =\m_1 \cup -\m_2$ (where the minus sign denotes the change of
orientation), endowed with a submersion $q:W\rightarrow \s \times
[0,1]$, with $\m_1=q^{-1}(\s \times \{ 0 \})$ and $\m_2=q^{-1}(\s
\times \{ 1 \})$ and $q|_{\m_1}=\pi_1$ and $q|_{\m_2}=\pi_2$. One
also requires that the map $q$ is a covering map branched along a
surface $S \subset \s\times [0,1]$ such that $\partial S =E_1 \cup
-E_2$, with $E_1=S\cap(\s \times \{ 0 \})$ and $E_2=S\cap (\s \times
\{ 1 \})$.
Since in the case of both $3$-manifolds and $4$-manifolds there is no
substantial difference in working in the $PL$ or smooth categories,
we keep formulating everything in the $PL$ setting.
We adapt easily this notion to the case of our correspondences. We
simply need to modify the definition above to take into account the
fact that our correspondences have two (not just one) covering maps
to $\s$, so that the cobordisms have to be chosen accordingly.
\begin{defn}\cite{AMM}\label{cobord}
Suppose given two morphisms $\m_1$ and $\m_2$ in $Hom(G,G')$, of
the form
$$ G \subset E_1 \subset \s \stackrel{\pi_{G,1}}{\longleftarrow} \m_1
\stackrel{\pi_{G',1}}{\longrightarrow}  \s \supset E'_1 \supset G'
$$
$$ G \subset E_2 \subset \s \stackrel{\pi_{G,2}}{\longleftarrow} \m_2
\stackrel{\pi_{G',2}}{\longrightarrow}  \s \supset E'_2 \supset G'.
$$ Then a cobordism between $\m_1$ and $\m_2$ is a $4$-dimensional
manifold $W$ with boundary $\partial W= \m_1 \cup -\m_2$, endowed
with two branched covering maps
\begin{equation}\label{Wmaps}
S \subset \s\times [0,1] \stackrel{q}{\longleftarrow} W
\stackrel{q'}{\longrightarrow} \s\times [0,1] \supset S',
\end{equation}
branched along surfaces $S,S'\subset \s\times [0,1]$. The maps $q$
and $q'$ have the properties that $\m_1 = q^{-1}(\s\times \{0\})=
q'^{-1}(\s\times \{0\})$ and $\m_2 =q^{-1}(\s\times \{1\})=
q'^{-1}(\s\times \{1\})$, with $q|_{\m_1}=\pi_{G,1}$,
$q'|_{\m_1}=\pi_{G',1}$, $q|_{\m_2}=\pi_{G,2}$ and
$q'|_{\m_2}=\pi_{G',2}$. The surfaces $S$ and $S'$ have boundary
$\partial S= E_1\cup -E_2$ and $\partial S'=E_1'\cup -E_2'$, with
$E_1= S \cap (\s\times \{0\})$, $E_2=S \cap (\s\times \{1\})$,
$E_1'=S' \cap (\s\times \{0\})$, and $E_2'=S' \cap (\s\times
\{1\})$.
\end{defn}
Here By ``surface" we mean a $2$-dimensional simplicial complex that is
PL-embedded in $\s \times [0,1]$, with boundary $\partial S \subset \s \times \{0,1\}$
given by 1-dimensional simplicial complexes, \ie embedded graphs.
\subsection{The Functoriality of Khovanov-Kauffman homology and Floer-Kauffman Homology}\label{cobo2}
In sections (\ref{KKH}) and (\ref{FKH}) we discussed the idea of associating a graph Khovanov-Kauffman Homology homomorphism and a graph Floer-Kauffman Homology to a cobordism between two graphs $G$ and $G'$. Khovanov homology is usually constructed only for knots and links and even for graphs in $\s$. We need a version of a homological invariant in 3-manifold of knots, links and graphs that is naturally defined for any 3-manifold. One such possibility would be Floer homology which is conjecturally related to Khovanov homology in the case of the 3-sphere. Kauffman in his work pointed to the possibility of associating a family of links to a graph inside 3-dimensional manifold $\m$ and this work help us in addition to the covering map
to prove the existence of the subcomplexes in 3-manifold $\m$ to study Floer homology for graphs inside $\m$.
\begin{thm}
Let  $$ G \subset \s \stackrel{\pi}{\longleftarrow} \m \stackrel{\pi^{-1}}{\longrightarrow}  \s \supset G'.$$
We need to find a map between the two Floer Homology groups $(M,\pi,\pi^{-1})$
$$FKH (G,S^3) \overset{\Upsilon_{(M,\pi,\pi^{-1})}}{\longrightarrow} FKH (G',S^3) $$
\end{thm}
\proof we need to use the generators and boundaries explicitly
 \[
\begin{xy}
\xymatrix{
 FC_*(G,S^3) \ar[d]^{\partial_*} \ar[r]^{ \Phi_1} &  \ar[d]^{\partial_*} FC_*(\pi^{-1}(G),M)  \\
 FC_{*-1}(G,S^3)  \ar[r]^{\Phi_1}& FC_{*-1}(\pi^{-1}(G),M) }
\end{xy}
\]

$$ FC_*(G) \overset{\Psi_2}{\longrightarrow} FC_*(G \cup G') \overset{\Psi_3}{\longrightarrow} FC_*(G')$$

\[
\begin{xy}
\xymatrix{
 FC_*(\pi^{-1}(G'),M ) \ar[d]^{\partial_*} \ar[r]^{ \Phi'_4} &  \ar[d]^{\partial_*} FC_*(G',S^3)  \\
 FC_{*-1}(\pi^{-1}(G'),M)  \ar[r]^{\Phi'_4}& FC_{*-1}(G',S^3) }
\end{xy}
\]

$$FKH (G,S^3) \overset{\Upsilon_{(M,\pi,\pi^{-1})}}{\longrightarrow} FKH (G',S^3) $$

 where $$\Upsilon_{(M,\pi,\pi^{-1})}= \Phi'_4 \circ \Psi_3 \circ \Psi_2 \circ \Phi_1 $$
Or\\
 \[
\begin{xy}
\xymatrix{
 \widehat{CFK}_*(T(G),S^3) \ar[d]^{\partial_*} \ar[r]^{ \Phi_1} &  \ar[d]^{\partial_*} \widehat{CFK}_*(\pi^{-1}(T(G)),M)  \\
 \widehat{CFK}_{*-1}(T(G),S^3)  \ar[r]^{\Phi_1}& \widehat{CFK}_{*-1}(\pi^{-1}(T(G)T(G)),M) }
\end{xy}
\]

$$ \widehat{CFK}_*(T(G)) \overset{\Psi_2}{\longrightarrow} \widehat{CFK}_*(T(G) \cup T(G')) \overset{\Psi_3}{\longrightarrow} \widehat{CFK}_*(T(G'))$$

\[
\begin{xy}
\xymatrix{
 \widehat{CFK}_*(\pi^{-1}(T(G')),M ) \ar[d]^{\partial_*} \ar[r]^{ \Phi'_4} &  \ar[d]^{\partial_*} \widehat{CFK}_*(T(G'),S^3)  \\
 \widehat{CFK}_{*-1}(\pi^{-1}(T(G')),M)  \ar[r]^{\Phi'_4}& \widehat{CFK}_{*-1}(T(G'),S^3) }
\end{xy}
\]

$$\widehat{HFK}(T(G),S^3) \overset{\Upsilon_{(M,\pi,\pi^{-1})}}{\longrightarrow} \widehat{HFK}(T(G'),S^3) $$

 where $$\Upsilon_{(M,\pi,\pi^{-1})}= \Phi'_4 \circ \Psi_3 \circ \Psi_2 \circ \Phi_1 $$

\diagram & & & & \hat{\m}=\m\times_{G'} \tilde{\m} \dlto^{P_1}
\drto_{P_2} & & &
\\
& & & \m \dlto_{\pi_{G}} \drto_{{\pi}_{G'}} & & \tilde{\m}
\dlto^{\tilde{\pi}_{G'}} \drto^{{\tilde{\pi}}_{G''}} & & \\
& & G\subset \s & & G'\subset \s & & G''\subset \s
\enddiagram

The fibered product $\hat{\m}$ is by definition a subset of the
product $\m\times \tilde\m$ defined as the preimage $\hat{\m}=
(\pi_{G'} \times \pi_{G''})^{-1} (\Delta (\s))$, where $\Delta (\s)$
is the diagonal embedding of $\s$ in $\s\times \s$. This defines a
smooth 3-dimensional submanifold of $\m\times \tilde\m$.

\[
\begin{xy}
\xymatrix{
 \widehat{CFK}_*(G,S^3) \ar[d]^{\partial_*} \ar[r]^{ \Phi_1} &  \ar[d]^{\partial_*} \widehat{CFK}_*(P_1\circ\pi^{-1}(G), \m\times \tilde\m)  \\
 \widehat{CFK}_{*-1}(G,S^3)  \ar[r]^{\Phi_1}& \widehat{CFK}_{*-1}(P_1\circ\pi^{-1}(G), \m\times \tilde\m) }
\end{xy}
\]

$$ \widehat{CFK}_*(G) \overset{\Psi_2}{\longrightarrow} \widehat{CFK}_*(G \cup G') \overset{\Psi_3}{\longrightarrow} \widehat{CFK}_*(G')$$

\[
\begin{xy}
\xymatrix{
 \widehat{CFK}_*(P_2 \circ \tilde{\pi}_{G''}, \m\times \tilde\m ) \ar[d]^{\partial_*} \ar[r]^{ \Phi'_4} &  \ar[d]^{\partial_*} \widehat{CFK}_*(G'',S^3)  \\
 \widehat{CFK}_{*-1}(P_2\circ\tilde{\pi}_{G''}, \m\times \tilde\m)  \ar[r]^{\Phi'_4}& \widehat{CFK}_{*-1}(G'',S^3) }
\end{xy}
\]

$$\widehat{HFK}(G,S^3) \overset{\Upsilon_{( \m\times \tilde\m,P_1\circ\pi,P_2\circ\tilde{\pi})}}{\longrightarrow} \widehat{HFK}(G',S^3) $$

 where $$\Upsilon_{(\m\times \tilde\m,P_1\circ\pi,P_2\circ\tilde{\pi})}= \Phi'_4 \circ \Psi_3 \circ \Psi_2 \circ \Phi_1 $$

 \endproof

 We need to use Hilden and Little definition in section (\ref{HL}) of a branched surface $\Sigma_P$ in $\s \times [0,1]$ by using a branched covering map $q: W \to \s \times [0,1]$, where $W = \m \times [0,1]$ a 4-dimensional Manifold, this branched map induces later two functors from the graph Floer-Kauffman Homology category $\mathcal{C}_{FKh}$ to the graph Khovanov-Kauffman Homology $\mathcal{C}_{KKh}$ category, and the second one from the Floer Homology category to the  Khovanov Homology category.

Consider $\m$ in $Hom(G,G')$ specified by a diagram $$ G \subset E \subset
\s \stackrel{\pi_{1}}{\longleftarrow} \m
\stackrel{\pi_{2}}{\longrightarrow}  \s \supset E' \supset G'.$$
 We can choose $W=\m \times [0,1]$ as a cobordism of $\m$ with itself. This has
$\partial W=\m \cup -\m$ with covering maps
$$ G \subset \s \times[0,1] \stackrel{q|_{\m \times \{0\}}}{\longleftarrow} W=\m \times[0,1]
\stackrel{q|_{\m \times \{1\}}}{\longrightarrow}  \s\times[0,1] \supset G'$$
branched along the PL surfaces $\Sigma_P$ in $\s\times [0,1]$ with $\partial \Sigma_P= G  \cup -G'$.
The branched covering maps $q|_{\m \times \{0\}}=\pi_{1}$ and $q|_{\m \times \{1\}}= \pi_{2}$ have the
properties that
$$ \m = q_1^{-1}(\s\times \{0\})=q_1^{-1}(\s\times \{1\})$$ Now we have $\pi_1^{-1}(G), \pi_2^{-1}(G')$ two graphs (not necessary connected) in $\m$.
Definition (\ref{branch}) showed that the inverse of image of the branched set is one-dimensional subcomplex and this
help us to define a new cobordism $\Sigma_M$ in a 3-manifold $\m$ in $W= \m \times [0,1]$ with $\partial \Sigma_M = \pi_1^{-1}(G) \cup -\pi_2^{-1}(G')$

  \[
\begin{xy}
\xymatrix{
\pi_1^{-1}(G)   \ar[d]^{\pi_1} \ar[r]^{\Sigma_M } &  \ar[d]^{\pi_2} \pi_2^{-1} (G')  \\
 G  \ar[r]^{\Sigma_P}& G' }
\end{xy}
\]
more precisely,
  \[
\begin{xy}
\xymatrix{
\pi_1^{-1}(G)  \subset \m \times \{0\} \ar[d]^{\pi_1} \ar[r]^{\m \times [0, 1]}_{\Sigma_M} &  \ar[d]^{\pi_2} M^3 \times \{1\} \supset \pi_2^{-1} (G')\\
 G \subset \s \times \{0\} \ar[r]^{\s \times [0, 1]}_{\Sigma_P}&  \s \times \{1\} \supset G' }
\end{xy}
\]

Now we can study the association of the Floer-Kauffman Homology for graphs in $\pi_1^{-1}(G)$ and $ \pi_2^{-1}(G')$ in a 3-dimensional manifold  $\m$
and study Khovanov-Kauffman Homology for graphs in $G, G'$ in a 3-dimensional manifold  $\s$
In $\m$ we found a graphs $\pi_1^{-1} (G), \pi_2^{-1} (G')$ under the inverse image of the branched cover maps $q|_{\m \times \{0\}}=\pi_{1}$ and $q|_{\m \times \{1\}}= \pi_{2}$ with $\Sigma_M$ cobordism relating these two graphs. In a parallel line $\Sigma_P$ is defined in $\s \times [0,1]$ as a branched surface by the map
$q: W= \m \times [0,1] \to \s \times [0,1] $ to be a cobordism with boundary $\partial \Sigma_P= G\cup -G'$. Our question is the relation between these two cobordisms $\Sigma_M$ and $\Sigma_P$, \ie it is possible to to use the branched cover map $q(\Sigma_M) \equiv \Sigma_P$?. These cobordisms induced Homomorphisms for both Khovanov-Kauffman homology for graphs $f_{KKh}: KKh(G) \to KKh(G')$ where $G$ and $G'$ are embedded graphs in $\s$ and one for Floer-Kauffman Homology for graphs $f_{FKh}: \widehat{HFG}(\pi_1^{-1} (G)) \to \widehat{HFG}(\pi_2^{-1} (G'))$ where $\pi_1^{-1}(G)$ and $ \pi_2^{-1}(G')$ are embedded graphs in a 3-manifold $\m$. The branched cover map $q$ induced a functor $\phi$ from the graph Floer-Kauffman Homology category $\mathcal{C}_{FKh}$ to the graph Khovanov-Kauffman Homology category $\mathcal{C}_{KKh}$, which takes the object $\widehat{HFG}((G))$ in $\mathcal{C}_{FKh}$ to the object  $KKh(G)$ in $\mathcal{C}_{KKh}$ and morphism $f_{FKh}$ to the morphism $f_{KKh}$ \ie $\phi(f_{FKh})= f_{KKh} $

        \[
\begin{xy}
\xymatrix{
FKh(\pi_1^{-1}(G))   \ar[d]^{\phi} \ar[r]^{ f_{FKh}}_{\Sigma_M} &  \ar[d]^{\phi} FKh(\pi_2^{-1}(G'))  \\
 KKh(G)  \ar[r]^{f_{KKh}}_{\Sigma_P}& KKh(G') }
\end{xy}
\]

If we use the idea of Kauffman of associating family of links to each graph and use the family of smooth cobordisms to the PL cobordism, and hence we can think by a another functor call it $\psi$ from the link Floer Homology category to the link Khovanov Homology category, which takes the object $\widehat{HFK}((L))$ for a link $L$ to the object  $Kh(L)$ and morphism $f_{Fh}$ to the morphism $f_{Kh}$
  \[
\begin{xy}
\xymatrix{
 T(\pi_1^{-1}(G)) \ar[d]^{\pi_1} \ar[r]^{ T(S_\alpha)_{\alpha \in \lambda}} &  \ar[d]^{\pi_2} T(\pi_2^{-1}(G'))  \\
 T(G)  \ar[r]^{T(\Sigma_\alpha)_{\alpha \in \lambda}}& T(G') }
\end{xy}
\]
Here we found a relation between Floer Homology and Khovanov Homology.

\bibliographystyle{amsplain}

\end{document}